\documentclass[fleqn,11pt]{article}
\newcommand{\scriptC}{$ \mathcal{C}\: $}
\newcommand{\boxtheta}{\{\theta\}}
\usepackage{graphicx}

\title{Axiomatic Characterization of Ordinary Differential Cohomology}
\author{James Simons and Dennis Sullivan}
\date{}

\begin{document}
\maketitle

\section*{Introduction}

The ring of differential characters, $\hat{H}(M,R/Z)$, a graded functor on the category
of smooth manifolds together with smooth maps, was developed in [6],
[7] and [8].  The motivation was to provide a home in
the base for the bundle invariants constructed in [5].

$\hat{H}$ satisfies the Character Diagram shown in \S1 and thus has natural
transformations onto both closed differential forms with integral periods and
integral cohomology.  It was shown in [6] that an enriched version of the Weil
homomorphism, carrying both an invariant polynomial and an associated integral
homology class of the classifying space, naturally factors through $\hat{H}$ on its way
to these targets.  Based on this and other considerations, $\hat{H}$ has been shown to
have useful applications to a number of areas, including conformal geometry and
flat bundles.  Most recently, the characters, perhaps generalized to relate to
K-theory, have appeared in descriptions of quantizations of field theory related
to super gravity and string theory, cf. [14].

In roughly the same time frame that differential characters were
introduced, Deligne cohomology was developed as a tool in algebraic
geometry, cf. [11].  When subsequently applied to the
smooth category it was shown to yield a
functor essentially equivalent to $\hat{H}$, cf. [11].  Thereafter other similar functors were developed, both in the smooth
category and others, cf. [9], [10], [12], [13].  Those in the smooth category all satisfy the Character Diagram, and,
while constructed in distinct manners, all have been shown to be naturally equivalent to $\hat{H}$.

Our purpose here is to show that, in the smooth category, the Diagram itself is
sufficient to uniquely characterize all such functors.  Theorem 1.1 shows that
$\hat{H}$ is unique up to a unique natural equivalence, and Theorem 1.2 shows that this
equivalence must preserve any graded product structure compatible with the
natural pairings on the adjacent functors of the Diagram.  The arguments are
based on naturality and on largely classical facts about perturbing arbitrary singular
cycles and homologies into embedded pseudomanifolds, both closed and
with boundary, cf. [1], [2], [3], [4].

In [14] an analog of differential characters, called differential
cohomology, is constructed for any generalized cohomology theory
(meaning the Eilenberg Steenrod axioms are satisfied except that the
groups of a point may not be concentrated in degree zero).  Each
flavor of differential cohomology satisfies an analog of the Character
Diagram, which, in the ordinary case, agrees with ours.  It would be
interesting to know whether each of the extraordinary flavors is
uniquely characterized by its diagram perhaps augmented by further
axioms, i.e. whether or not there is a version of Theorem 1.1 for
generalized differential cohomology.

We are very happy to thank Jeff Cheeger for useful conversations in the early
stages of the work, Blaine Lawson for acquainting us with the current status
of the field, and Chris Bishop for his help with the proof of Fact 2.1.

\section*{\S1. Uniqueness Theorems}

The proofs presented in this section depend on results outlined in \S2.

Let \scriptC denote the category of smooth manifolds together with
smooth maps. Let $\Lambda$
denote the functor from \scriptC to differential graded algebras given
by smooth real valued differential forms with the usual $d$ operator and
wedge product.  Let $\Lambda_{Z} \subseteq \Lambda$ denote the
functor which assigns to any manifold the graded algebra
of closed differential forms with integral periods, and
$\Lambda/\Lambda_{Z}$ be the quotient functor, taking values in graded
abelian groups.  Finally, let
$H(A)$ denote the smooth singular cohomology functor with coefficients
in $A$, an abelian group or commutative ring.

\underline{\textbf{Definition}}: A \textbf{Character Functor} is a
5-tuple $\{\hat{G},i_{1},i_{2},\delta_{1},\delta_{2}\}$, where
$\hat{G}$ is a functor from \scriptC to graded abelian groups, and
$i_{1},i_{2},\delta_{1},\delta_{2}$ are natural transformations which
make the following \textbf{Character Diagram} commutative and its
diagonal sequences exact.

\vspace{.5cm}
\begin{center}
\setlength{\unitlength}{0.5cm}
\begin{picture}(24,16)\thicklines
\put(4,0){$0$}
\put(21,0){$0$}
\put(5.5,1){\vector(1,1){1.5}}
\put(18.5,3){\vector(1,-1){1.5}}
\put(8,4){$\Lambda^{k-1}/\Lambda^{k-1}_{Z}$}
\put(12.5,4){\vector(1,0){2.5}}
\put(17,4){$\Lambda^{k}_{Z}$}
\put(13.5,4.5){\small{$d$}}
\put(2.5,6){\vector(1,1){1.5}}
\put(6.5,7){\vector(1,-1){1.5}}
\put(7.5,7){\small{$\beta$}}
\put(10.5,7){\small{$i_{2}$}}
\put(10.5,6){\vector(1,1){1.5}}
\put(14.5,7){\vector(1,-1){1.5}}
\put(15.0,7){\small{$\delta_{1}$}}
\put(18.5,7){\small{$s$}}
\put(18.5,6){\vector(1,1){1.5}}
\put(22.5,7){\vector(1,-1){1.5}}
\put(4,8){$H^{k-1}(R)$}
\put(13,8){$\hat{G^{k}}$}
\put(20,8){$H^{k}(R)$}
\put(2.5,11){\vector(1,-1){1.5}}
\put(6.5,11){\small{$\alpha$}}
\put(6.5,10){\vector(1,1){1.5}}
\put(10.5,11){\vector(1,-1){1.5}}
\put(11.5,10.5){\small{$i_{1}$}}
\put(14,10.5){\small{$\delta_{2}$}}
\put(14.0,9.5){\vector(1,1){1.5}}
\put(18.5,11){\vector(1,-1){1.5}}
\put(19.5,10.5){\small{$r$}}
\put(22.5,10){\vector(1,1){1.5}}
\put(8,12){$H^{k-1}(R/Z)$}
\put(13,12){\vector(1,0){2}}
\put(16,12){$H^{k}(Z)$}
\put(13.5,12.5){\small{$-B$}}

\put(5.5,15){\vector(1,-1){1.5}}
\put(18.5,14){\vector(1,1){1.5}}
\put(4,16){$0$}
\put(21,16){$0$}
\end{picture}
\end{center}

where, $(\alpha, B, r)$ is the Bockstein long exact sequence
associated to the coefficient sequence
\mbox{$Z \rightarrow R \rightarrow R/Z$}, and $(\beta, d, s)$ is another long exact sequence in which $\beta$ and
$s$ are obviously defined via the de Rham Theorem.

\underline{\textbf{Remark}}:  It will be shown in the proof of Theorem
1.1 and in Corollary 1.1 that
$\delta_{2}$ is actually uniquely determined by the other three
natural transformations, and thus, strictly speaking, is not required
in the above definition.  Nonetheless, its inclusion completes the picture.

\underline{\textbf{Differential Characters}}:  As mentioned in the
Introduction, differential characters $\hat{H}(M,R/Z)$ and their
related mappings comprise a character functor.  We recall their
construction:

For $M \in \mathcal{C}$ let $C_{k}(M)$ and $Z_{k}(M)$ respectively
denote the groups of smooth singular $k$-chains and $k$-cycles.  For
$\omega \in \Lambda^{k}(M)$ and $a \in C_{k}(M)$ we define $\tilde{\omega} \in C^{k}(M,R/Z)$ by 

\begin{displaymath}
1.1) \qquad \tilde{\omega}(a) = \int_{a}\omega \; \bmod Z.
\end{displaymath}

Following
\renewcommand{\thefootnote}{\fnsymbol{footnote}}
[8]\footnote[1]{The indexing has been shifted by 1 from that in [8].} 
\renewcommand{\thefootnote}{\arabic{footnote}}
we define

\begin{displaymath}
1.2) \qquad \hat{H}^k(M,R/Z) = \{ f \in \textrm{Hom}(Z_{k-1}(M),R/Z)
\quad|\quad f \circ \partial = \tilde{\omega}_{f} \}
\end{displaymath}

for some $\omega_{f} \in \Lambda^{k}(M)$.  It is easily seen that
$\omega_{f}$ is uniquely determined by $f$ and that $\omega_{f} \in
\Lambda^{k}_{Z}(M)$.  In the notation of the Character Diagram

\begin{displaymath}
1.3) \qquad \delta_{1}(f) = \omega_{f}.
\end{displaymath}

Since $H^{k-1}(M,R/Z) \cong \textrm{Hom}(H_{k}(M),R/Z)$ we may
consider $H^{k-1}(M,R/Z) \subseteq \hat{H}^{k}(M,R/Z)$.  From 1.1) and
Stokes Theorem we see that $\theta \rightarrow \tilde{\theta}$ maps
$\Lambda^{k-1}(M) \rightarrow \hat{H}^{k}(M,R/Z)$, and since the
kernel of this map is $\Lambda^{k-1}_{Z}(M)$ we may consider
$\Lambda^{k-1}(M)/\Lambda^{k-1}_{Z}(M) \subseteq \hat{H}^{k}(M,R/Z)$.

\begin{displaymath}
1.4) \qquad \textrm{The above inclusions define } i_{1} \textrm{ and } i_{2}.
\end{displaymath}

Since $R$ is divisible we can find $T \in C^{k-1}(M,R)$ with
$\tilde{T}|Z_{k-1}(M) = f$.  It is easily shown that $\omega_{f} -
\delta T \in C^{k}(M,Z)$ and is closed.  Moreover its cohomology
class, $\{\omega_{f} - \delta T\} \in H^{k}(M,Z)$, is independent of
the choice of $T$.  By definition, $\delta_{2}(f) = \{\omega_{f} - \delta T\}$.

If $\phi: M \rightarrow N$ is $C^{\infty}$, $a \in Z_{k-1}(M)$ and $f
\in \hat{H}^{k}(M,R/Z)$ then $\phi^{*} : \hat{H}^{k}(N,R/Z)
\rightarrow \hat{H}^{k}(M,R/Z)$ is simply defined by

\begin{displaymath}
1.5) \qquad \phi^{*}(f)(a) = f(\phi_{*}(a)).
\end{displaymath}

\pagebreak

\underline{\textbf{Theorem 1.1}}: Any character functor $\{\hat{G},
i_{1}, i_{2}, \delta_{1}, \delta_{2}\}$ is equivalent to $\{\hat{H},
i_{1}, i_{2}, \delta_{1}, \delta_{2}\}$ via a natural transformation
$\Phi : \hat{G} \rightarrow \hat{H}$ which commutes with the identity
map on all other functors in the Diagram; such $\Phi$ is unique.

\underline{\textbf{Proof}}:  

For any singular chain, $x$, we will denote by $|x|$ the union of the
images of all the simplexes of $x$.

Let $a \in Z_{k-1}(M)$.  Using Fact 2.1 we may choose a $(k-1)$-good
open neighborhood $U$ of $|a|$.  Let $\lambda : U \rightarrow M$
denote the inclusion map.  Since $H^{k}(U,Z) = 0$, the Character Diagram shows

\begin{displaymath}
1.6) \qquad \lambda^{*}(g) = i_{2}( \boxtheta )
\end{displaymath}

where $\theta \in \Lambda^{k-1}(U)$ and is determined up to an element
of $\Lambda^{k-1}_{Z}(U)$.  We set

\begin{displaymath}
1.7) \qquad \Phi(g)(a) = \tilde{\theta}(a) = \int_{a}\theta \; \bmod Z.
\end{displaymath}

We proceed by showing

\[
\begin{array}{rl}
1.8) & \Phi \textrm{ is a well defined homomorphism from }
\hat{G}^{k}(M) \textrm{ into } \textrm{Hom}(Z_{k}(M),R/Z). \\
\\
1.9) & \textrm{Im}(\Phi) \subseteq  \hat{H}^{k}(M,R/Z) \subseteq
\textrm{Hom} (Z_{k}(M),R/Z). \\
\\
1.10) & \textrm{If } \varphi:M \rightarrow N \textrm{ is } C^{\infty},
\textrm{ then } \Phi \circ \varphi^{*} = \varphi^{*} \circ \Phi,
\textrm{where } \varphi^{*} \textrm{ denotes the maps } \\
 & \textrm{induced by } \varphi
\textrm{ of the functors } \hat{G} \textrm{ and } \hat{H}. \\
\\
1.11) & \Phi \circ i_{1} = i_{1}, \: \Phi \circ i_{2} = i_{2}, \: \delta_{1}
\circ \Phi = \delta_{1}. \\
\\
1.12) & \Phi \textrm{ is an isomorphism.} \\
\\
1.13) & \delta_{2} \circ \Phi = \delta_{2}. \\
\\
1.14) & \Phi \textrm{ is the unique natural transformation from }
\hat{G} \textrm{ to } \hat{H} \textrm{ satisfying } \\
 & \textrm{1.11), 1,12) and 1.13).} \\
\end{array}
\]

To show 1.8), note that clearly 1.7) is independent of the choice of $\theta \in \boxtheta$.
To show it is independent of the choice of $U$, let $U^{\prime}$ be
another $(k-1)$-good neighborhood of $|a|$.  Again by Fact 2.1 we can
choose a third $(k-1)$-good neighborhood of $|a|$, $U^{\prime\prime}
\subseteq U \cap U^{\prime}$.  Let $U^{\prime\prime} \stackrel{\eta}{\longrightarrow} U
\stackrel{\lambda}{\longrightarrow} M$ be the inclusion maps.  By 1.5) and the naturality of $\hat{G}$

\begin{displaymath}
(\lambda \circ \eta)^{*}(g) = \eta^{*}(\lambda^{*}(g)) = \eta^{*}(i_{2}(\{\theta\}))
 = i_{2}(\eta^{*}(\{\theta\})) = i_{2}(\{\eta^{*}(\theta)\}).
\end{displaymath}

Since $\widetilde{\eta^{*}(\theta)}(a) = \tilde{\theta}(a)$ we
see that the definitions of $\Phi$ using $U$ and $U^{\prime\prime}$
agree.  Since the same is true for $U^{\prime}$ and $U^{\prime\prime}$, we have shown that the
definition of $\Phi$ is independent of the choice of $U$.

To show $\Phi(g) \in \textrm{Hom}(Z_{k-1}, R/Z)$ we must show that for $a,b
\in Z_{k-1}$ that $\Phi(g)(a+b) = \Phi(g)(a) + \Phi(g)(b)$.  Using all
of the distinct individual simplices of $a$ and $b$, create a chain
$c$ such that $|a| \cup |b| \subseteq |c|$.  Pick $U$, a $(k-1)$-good
neighborhood of $|c|$.  Since $|a+b| \subseteq |a| \cup |b|$ we have
each of $|a|$, $|b|$, $|a+b| \subseteq U$.  Choosing $\theta$ as in 1.6)

\begin{displaymath}
\Phi(g)(a+b) = \tilde{\theta}(a+b) = \tilde{\theta}(a) +
\tilde{\theta}(b) = \Phi(g)(a) + \Phi(g)(b).
\end{displaymath}

Since it is clear that $\Phi(g_{1} + g_{2}) = \Phi(g_{1}) +
\Phi(g_{2})$ we have now shown 1.8).  

To prove 1.9) we will show that if $a \equiv \partial e$ then

\begin{displaymath}
*)\qquad \Phi(g)(a) = \widetilde{\delta_{1}(g)}(e).
\end{displaymath}

Since both sides vanish if $k-1 = \textrm{ dim } M$, we may assume
$k-1 < \textrm{ dim } M$.  Let $U$ be a $(k-1)$-good neighborhood of
$|a|$.  Using Fact 2.2 we can find a $(k-1)$ dim imbedded
pseudomanifold $P \subseteq U$ and a $k$-chain $b$ with $|b| \subseteq U$ such that

\begin{displaymath}
a = \partial b + P
\end{displaymath}

where we identify $P$ with its fundamental cycle.  Letting $\theta$
be defined as in 1.6) we see
\begin{eqnarray*}
1.15)\quad \Phi(g)(a) & = & \tilde{\theta}(\partial b) + \Phi(g)(P) =
\widetilde{d\theta}(b) + \Phi(g)(P) = \widetilde{\delta_{1}(\lambda^{*}(g))}(b) +
\Phi(g)(P) \\
& = & \widetilde{\delta_{1}(g)}(b) + \Phi(g)(P)
\end{eqnarray*}
where we used Stokes Theorem, the fact that $d = \delta_{1} \circ
i_{2}$, and naturality of $\hat{G}$.

Since $a$ is a boundary in $M$ and $P$ is homologous to $a$, $P$ is a
boundary in $M$.  We then use Fact 2.3 to find $U^{\prime}$, a $(k-1)$-good neighborhood of $P$ with $\lambda : U^{\prime} \rightarrow N$ the
inclusion, and a $k$-chain $y$ with $ |y| \subseteq U^{\prime}$ and $P = \partial y$.

As in 1.6), $\lambda^{*}(g) = i_{2}\{\theta^{\prime}\} $ for
some $\theta^{\prime} \in \Lambda^{k-1}(U^{\prime})$, and by the same
argument as above

\begin{displaymath}
\Phi(g)(P) = \tilde{\theta}^{\prime}(\partial y) = \widetilde{\delta_{1}(g)}(y).
\end{displaymath}

Combining this with 1.15) we see

\begin{displaymath}
\Phi(g)(a) = \widetilde{\delta_{1}(g)}(b+y) = \widetilde{\delta_{1}(g)}(b
+y-e+e) = \widetilde{\delta_{1}(g)}(b+y-e) + \widetilde{\delta_{1}(g)}(e).
\end{displaymath}

But, since $\partial(b+y)=a=\partial e$, $b+y-e$ is a cycle, and
$\widetilde{\delta_{1}(g)} \in \Lambda^{k}_{Z}(M)$, it follows the first
term vanishes.  This proves $\ast$) and thus 1.9).

To prove 1.10), let $g \in \hat{G}^{k}(N)$, $a \in Z_{k-1}(M)$, $U
\subseteq N$, a $(k-1)$-good neighborhood of $|\phi_{\ast}(a)|$ and $W
\subseteq \phi^{-1}(U)$, a $(k-1)$-good neighborhood of $|a|$.  The
various inclusions are labelled in the commutative diagram

\[
\begin{array}{ccccc}
W & \stackrel{\gamma}{\subset} & \phi^{-1}(U) & \stackrel{\eta}{\subset} &
M \\
& & \downarrow \phi & & \downarrow \phi \\
& & U & \stackrel{\lambda}{\subset} & N \\
\end{array}
\]

Now, $\lambda^{\ast}(g) = i_{2}(\{\theta\})$ for $\theta \in
\Lambda^{k-1}(U)$.  By naturality of $g$

\begin{displaymath}
(\eta \circ \gamma)^{\ast}(\phi^{\ast}(g)) = \gamma^{\ast} \circ
\eta^{\ast} \circ \phi^{\ast}(g) = \gamma^{\ast} \circ \phi^{\ast}
\circ \lambda^{\ast}(g) = (\phi \circ
\gamma)^{\ast}(i_{2}(\{\theta\})) = i_2(\{(\phi \circ \gamma)^{\ast}(\theta)\}).
\end{displaymath}

Thus, by the definition of $\Phi$ and 1.5)

\begin{displaymath}
(\Phi(\phi^{\ast}(g)))(a) = \widetilde{(\phi \circ
\gamma)^{\ast}(\theta)}(a) = \tilde{\theta}((\phi \circ
\gamma)_{\ast}(a)) = \tilde{\theta}(\phi_{\ast}(a)) =
(\Phi(g))(\phi_{\ast}(a)) = \phi^{\ast}(\Phi(g))(a).
\end{displaymath}

To prove 1.11) let $\mu \in H^{k-1}(M,R/Z)$, $a \in Z_{k-1}$ and $U$ a
$(k-1)$-good neighborhood of $|a|$.  Since $H^{k}(U) = 0$, by the
exactness of the Bockstein sequence, $\mu = \alpha(x)$, where $x \in
H^{k-1}(R)$.  By commutativity of the Character Diagram for
$\hat{G}$

\begin{displaymath}
i_{i}(\mu) = i_{i}(\alpha(x)) = i_{2}(\beta(x)).
\end{displaymath}

Thus, by 1.6)

\begin{displaymath}
\Phi(i_{1}(\mu))(a) = \tilde{\theta}(a) \qquad \textrm{ for any }
\theta \in \{{\beta(x)}\}.
\end{displaymath}

But, since $\beta$ is defined by the de Rham Theorem, and using $i_{1}$
of $\hat{H}$ as defined in 1.4)

\begin{displaymath}
\tilde{\theta}(a) = x(a) \bmod Z = \alpha(x)(a) = \mu(a) = i_{1}(\mu)(a)
\end{displaymath}

Thus $\Phi \circ i_{1} = i_{1}$.  That $\Phi \circ i_{2} = i_{2}$ follows immediately from 1.5) and
1.6), and that $\delta_{1} \circ \Phi = \delta_{1}$ follows from $\ast$) in
the proof of 1.8).  This shows 1.11).

To prove 1.12) we note that the following diagram of exact sequences is commutative:

\[
\begin{array}{ccccccccc}
0 & \longrightarrow & H^{k-1}(M,R/Z) & \stackrel{i_{1}}{\longrightarrow} &
\hat{G}^{k}(M) & \stackrel{\delta_{1}}{\longrightarrow} & \Lambda^{k}_{Z}(M) &
\longrightarrow & 0 \\
& & {\scriptstyle \|} & & \downarrow \Phi & & {\scriptstyle\|} & & \\
0 & \longrightarrow & H^{k-1}(M,R/Z) & \stackrel{i_{1}}{\longrightarrow} &
H^{k}(M,R/Z) & \stackrel{\delta_{1}}{\longrightarrow} & \Lambda^{k}_{Z}(M) &
\longrightarrow & 0 \\
\end{array}
\]

Thus by the Five Lemma $\Phi$ is an isomorphism.

To prove 1.13) we need

\underline{\textbf{Lemma 1.1}}:  For any character functor, $\hat{G}$,
$\delta_{2}$ is uniquely determined by $i_{1}, i_{2}$ and $\delta_{1}$.

\underline{\textbf{Proof}}:  Let $\delta^{\prime}_{2}:\hat{G}^{k}(Z)
\rightarrow H^{k}(Z)$
be another such $\delta_{2}$, which, together with $i_{1}, i_{2}$ and
$\delta_{1}$ make $\hat{G}$ a character functor.  Since each of
$\delta_{2}$ and $\delta^{\prime}_{2}$ induce isomorphisms from
$\hat{G}^{k}/i_{2}(\Lambda^{k-1}/\Lambda^{k-1}_{Z})$ onto $H^{k}(Z)$,

\begin{displaymath}
\delta^{\prime}_{2} \circ \delta^{-1}_{2} : H^{k}(Z) \rightarrow H^{k}(Z)
\end{displaymath}

is an automorphism.

Let $T^{k}(Z) \subseteq H^{k}(Z)$ denote the torsion subgroup.  For
$\tau \in T^{k}(Z)$, the Bockstein exact sequence shows $\tau =
B(\mu)$ for some $\mu \in H^{k-1}(R/Z)$.  Thus, by the commutativity
of the Character Diagram

\begin{displaymath}
\delta^{\prime}_{2} \circ \delta^{-1}_{2}(\tau) = \delta^{\prime}_{2}
\circ \delta^{-1}_{2}(B(\mu)) = \delta^{\prime}_{2} \circ
\delta^{-1}_{2}(\delta_{2} \circ i_{1}(\mu)) = \delta^{\prime}_{2}
\circ i(\mu) = B(\mu) = \tau.
\end{displaymath}

Therefore $(\delta^{\prime}_{2} \circ \delta^{-1}_{2})|T^{k}(Z)$ is
the identity.

By naturality of $\hat{G}$, equipped with either $\delta_{2}$ or
$\delta_{2}^{\prime}$, we see that $\delta^{\prime}_{2} \circ
\delta^{-1}_{2}$ is a natural automorphism of $H^{k}(Z)$ which holds
fixed $T^{k}(Z)$.  By Fact 1.1 below this can only hold if
$\delta_{2}^{\prime} \circ \delta_{2}^{-1}$ is the identity.  Thus
$\delta_{2}^{\prime} = \delta_{2}$. $\diamondsuit$

\underline{\textbf{Fact 1.1}}:  Any natural (homotopy invariant)
automorphism of the integral cohomology functor on \scriptC which is
either identity on torsion, or the identity after tensoring with $R$
must be the identity.

\underline{\textbf{Proof}}:  

One knows any integral $k$-cohomology class of a finite
dimensional manifold $M$ is induced by pulling back a universal class 
$u(k)$ on a fixed  universal space $K(Z,k)$, by a map  $f$ of $M$ into a
finite skeleton of $K(Z,k)$.  The space $K(Z,k)$ is characterized by having
one non zero homotopy group $Z$ in dimension $k$.  This finite skeleton
of $K(Z,k)$ has the same homotopy type as a manifold $V$, as can be
seen by embedding the skeleton into a euclidean space and forming a
regular neighborhood.  The map into this smooth neighborhood can be
deformed to a smooth map which we also denote by $f$.  Now we
calculate the automorphism.  The $k$-th cohomology of this
manifold  neighborhood $V$ is $Z$ in dimension $k$, by the Hurewicz and
universal coefficient theorems applied to $K(Z,k)$.  It follows any
natural automorphism must multiply the universal class $u(k)$ and
thus the general  class  $f^{*}(u(k))$ by plus or minus one.
The minus one possibility is ruled out by either one of our further
hypotheses by considering a manifold with a nonzero odd order class
or an infinite order class in degree $k$. $\diamondsuit$

To prove 1.13) we set $\delta_{2}^{\prime} = \delta_{2} \circ \Phi$
and note that using 1.10), 1.11) and 1.12) one easily shows $\{\hat{G},
i_{1}, i_{2}, \delta_{1}, \delta_{2}^{\prime}\}$ is a character
functor.  By Lemma 1.1 we see $\delta_2 \circ \Phi = \delta_{2}$. 

To complete the proof of the theorem we must show 1.14).  If $\Phi^{\prime}$ were another such then $\Phi^{\prime}
\circ \Phi^{-1} : \hat{H}^{k}(M,R/Z) \rightarrow \hat{H}^{k-1}(M,R/Z)$
would be a natural automorphism holding fixed the other terms in the
Character Diagram.  Thus, for $f \in \hat{H}^{k-1}(M,R/Z)$, $a \in
Z_{k-1}(M)$, and $U$ a $(k-1)$-good open neighborhood of $|a|$ with
$\lambda : U \rightarrow M$,
\begin{eqnarray*}
(\Phi^{\prime} \circ \Phi^{-1})(f)(a) & = & (\Phi^{\prime} \circ
\Phi^{-1})(\lambda^{*}(f))(a) \\
 & = & (\Phi^{\prime} \circ \Phi^{-1})(i_{2}(\{\theta\})(a)) \\
 & = & i_{2}(\Phi^{\prime} \circ \Phi^{-1}(\{\theta\})(a)  \\
 & = & i_{2}(\{\theta\})(a) = (\lambda^{*}(f))(a) = f(a).  \\
\end{eqnarray*}
Thus $\Phi = \Phi^{\prime}$. $\diamondsuit$

As was indicated in the Remark following the definition of the
character functor, the result below shows that the definition actually
requires less data.

\underline{\textbf{Corollary 1.1}}

If $\hat{G}$ is a functor from \scriptC to graded abelian groups and
$i_{1}$, $i_{2}$ and $\delta_{1}$ are natural transformations making
the diagram below commute and its diagonal exact

\vspace{.5cm}
\begin{center}
\setlength{\unitlength}{0.5cm}
\begin{picture}(24,16)\thicklines
\put(21,0){$0$}
\put(18.5,3){\vector(1,-1){1.5}}
\put(8,4){$\Lambda^{k-1}/\Lambda^{k-1}_{Z}$}
\put(12.5,4){\vector(1,0){2.5}}
\put(17,4){$\Lambda^{k}_{Z}$}
\put(13.5,4.5){\small{$d$}}
\put(6.5,7){\vector(1,-1){1.5}}
\put(7.5,7){\small{$\beta$}}
\put(10.5,7){\small{$i_{2}$}}
\put(10.5,6){\vector(1,1){1.5}}
\put(14.5,7){\vector(1,-1){1.5}}
\put(15.0,7){\small{$\delta_{1}$}}
\put(4,8){$H^{k-1}(R)$}
\put(13,8){$\hat{G^{k}}$}
\put(6.5,11){\small{$\alpha$}}
\put(6.5,10){\vector(1,1){1.5}}
\put(10.5,11){\vector(1,-1){1.5}}
\put(11.5,10.5){\small{$i_{1}$}}
\put(8,12){$H^{k-1}(R/Z)$}
\put(5.5,15){\vector(1,-1){1.5}}
\put(4,16){$0$}
\end{picture}
\end{center}

then $i_{2}$ is $1:1$, and there exists one and only one
$\delta_{2}:\hat{G} \rightarrow H(Z)$ so that
$\{\hat{G},i_{1},i_{2},\delta_{1},\delta_{2}\}$ is a character
functor.

\underline{\textbf{Proof:}}  

$ \qquad i_{2}(\{\theta\}) = 0 \Rightarrow
\delta_{1} \circ i_{2}(\{\theta\}) = 0 \Rightarrow d\theta = 0
\Rightarrow \{\theta\} = \beta(x) \Rightarrow i_{1} \circ \alpha(x) =
0 \\
\Rightarrow \alpha(x) = 0 \Rightarrow x = r(u)$ for $u \in
H^{k-1}(Z)$, and thus $0 = \beta(x) = \{\theta\}$.  Thus $i_{2}$ is $1:1$.

We next show that if $U \subseteq M$ is $k$-good and $\lambda : U
\rightarrow M$ is the inclusion then, for $g \in \hat{G}^{k}(M)$

\begin{displaymath}
*)\qquad \lambda^{*}(g) = i_{2}(\{\theta\}) \qquad \textrm{ for unique
} \{\theta\}.
\end{displaymath}

To show $\ast$), note that since $U$ is $k$-good
$\delta_{1}(\lambda^{*}(g)) = d(\{\rho\})$ and thus
$\delta_{1}(i_{2}(\{\rho\})) = \delta_{1}(\lambda^{*}(g)) \Rightarrow
\lambda^{*}(g) = i_{2}\{\rho\} + i_{1}(u)$, where $u \in
H^{k}(U,R/Z)$.  Since $H^{k}(U,Z) = 0$ the Bockstein shows $u =
\alpha(x)$.  Thus $i_{1}(u) = i_{1}(\alpha(x)) = i_{2}(\beta(x))$.
Therefore $\lambda^{*}(g) = i_{2}(\{\rho\} + \alpha(x)) =
i_{2}(\{\theta\})$ for some $\{\theta\}$.  Since $i_{2}$ is $1:1$,
$\{\theta\}$ is unique.

We may then define $\Phi : \hat{G} \rightarrow \hat{H}$ as in 1.7),
and follow the proof of Theorem 1.1 through 1.12).  By then setting
$\delta_{2} = \delta_{2} \circ \Phi$, $\hat{G}$ becomes a character
functor, and the uniqueness of $\delta_{2}$ is assured by Lemma 1.1.
$\diamondsuit$

\pagebreak

\underline{\textbf{Ring Structure}}

In [7] $\hat{H}$ was shown to possess a natural associative graded
commutative ring
structure, $\ast$ , cf. [9]-[13].  For $f \in \hat{H}^{k}$ and $g \in \hat{H}^{l}$
this ring structure also satisfies, cf. [8],

\[
\begin{array}{rl}
1.16) & f \ast g \in \hat{H}^{k+l} \\
\\
1.17) & f \ast g = (-1)^{kl} \: g \ast f \\
\\
1.18) &  \delta_{1}(f \ast g) = \delta_{1}(f) \wedge \delta_{1}(g)  \\
\\
1.19) & \delta_{2}(f \ast g) = \delta_{2}(f) \cup \delta_{2}(g) \\
\\
1.20) & f \ast i_{1}(u) = (-1)^{k} \: i_{1} \: (\delta_{2}(f) \cup u) \\
\\
1.21) & f \ast i_{2}(\{\theta\}) = (-1)^{k} \: i_{2} \: (\{\delta_{1}(f) \wedge \theta \}) \\
\\

\end{array}
\]

\underline{\textbf{Theorem 1.2}}:  A character functor, $\hat{G}$,
possesses at most one natural ring structure satisfying 1.16) - 1.21).

\underline{\textbf{Proof}}:  Suppose $\ast$ and $\dag$ are two such.
Set $B(f,g) = f \ast g - f \, \dag \, g$.  From 1.16), 1.18) and the
Character Diagram we see that

\begin{displaymath}
B(f,g) \in \textrm{Ker} (\delta_{1}) = i_{1}(H^{k+l-1}(R/Z)).
\end{displaymath}

From 1.17) and 1.21) we see that if $f \in
i_{2}(\Lambda^{k-1}/\Lambda^{k-1}_{Z})$ or $g \in i_{2}(\Lambda^{l-1}/\Lambda^{l-1}_{Z})$

\begin{displaymath}
B(f,g) = 0.
\end{displaymath}

Since the diagram shows
$\hat{G}^{k}/i_{2}(\Lambda^{k-1}/\Lambda^{k-1}_{Z}) \cong H^{k}(Z)$ we
see

\begin{displaymath}
B: H^{k}(Z) \times H^{l}(Z) \rightarrow H^{k + l - 1}(R/Z).
\end{displaymath}

By the naturality assumptions on $\ast$ and $\dag$, we see that $B$ is
a natural transformation satisfying the hypotheses of Fact 1.2 below.  Thus
$B \equiv 0$ and $\ast = \dag$. $\diamondsuit$

\underline{\textbf{Fact 1.2}}:  Any natural (homotopy
invariant) cohomology operation on \scriptC which assigns to a pair of
integral cohomology classes $(a,b)$ in dimensions $(i,j)$ an $R/Z$
cohomology class $(a \cdot b)$ in dimension $i + j - 1$, so that $(0
\cdot b) = (a \cdot 0) = 0$ must be identically zero.

\pagebreak

\underline{\textbf{Proof}}:

We will tacitly replace finite skeleta of the spaces below that are
  sufficient for the calculation by manifold thickenings as we did in the proof of Fact 1.1.  Apply the  hypothetical operation to the universal classes
mentioned in the proof of Fact 1.1.  $u(k)$ and $u(j)$ pulled back by
the respective projections to the cartesian product $K(Z,k) \times K(Z,j)$.  We obtain a class $p$ in the product $K(Z,k) \times K(Z,j)$ in degree $k+j-1$ with
coefficients in $R/Z$.  To compute the hypothetical operation on a pair
of classes we use the pair of classes to build a map into $K(Z,k) \times K(Z,j)$ and
then pull back $p$.  By our hypothesis that $a \cdot 0 = 0$ we get that $p$
restricted to the second factor of $K(Z,k) \times K(Z,j)$ is zero.
  Similarly we deduce $p$ restricted to the first factor is zero.  It follows from the exact
sequence of a pair that $p$ comes from the smash product obtained by
collapsing the two axes in $K(Z,k) \times K(Z,j)$ to a point.
The smash product of a $j-1$ connected space and a $k-1$ connected space
is $k+j-1$ connected. So the $j+k-1$ cohomology with any coefficients  of
the smash product is zero.  Therefore $p$ is zero and Fact 1.2 is
  proved.  $\diamondsuit$

By Theorem 1.1 we may use $\Phi$ to pull back $\ast$ to obtain a ring
structure on $\hat{G}$.  We may thus combine Theorems 1.1 and 1.2 as

\underline{\textbf{Theorem 1.3}}:  Any character functor admits
exactly one ring structure satisfying 1.16) - 1.21).  Equipped with
this ring structure, any two character functors are naturally
equivalent by a unique equivalence consistent with the Character Diagram.

\underline{\textbf{Characters over $R/\Gamma$}}

Character Functors, together with a unique ring structure, may be
defined over $R/\Gamma$, where $\Gamma$ is any proper subring of $R$.
All the above definitions and results are true in this more general
setting, and the proofs are identical.

\section*{\S2 - Geometric Homology Background}

\underline{\textbf{Fact 2.1}}:  Let $K$ in $M$ denote the compact
image of a smooth singular $k$-chain in $M$.  Then every neighborhood
of $K$ contains a smaller neighborhood whose integral cohomology
vanishes above $k$.  (We call these $\mathbf{k}$\textbf{-good} neighborhoods.)

\underline{\textbf{Proof of 2.1}}:
\begin{itemize}
\item[i)] Choose a Riemannian metric on $M$.  Since maps satisfying the
Lipschitz condition, 

\begin{displaymath}
\textrm{distance } (f(x),f(y)) \leq L \cdot \textrm{distance}(x,y), 
\end{displaymath}

do not increase Hausdorf dimension [2], we see the Hausdorf dimension of $K$ is at most $k$.  In particular the
Hausdorf $(k+1)$ measure of $K$ is zero.

\item[ii)]  Since the Hausdorf  $k+1$ measure of $K$ is zero any local orthogonal projection of a compact piece of $K$ into a $k+1$
plane is a closed set of measure zero and so misses an open dense
set. Thus we can choose a small triangulation so that $K$ misses the
$n-k-1$ skeleton.  We do this by first choosing one small triangulation and then performing local translations
in $k+1$ directions orthogonal to the top cells of the $n-k-1$ skeleton into the open dense sets that miss the local
projections of $K$.  Then $K$ is contained in the handlebody
neighborhood of the dual $k$ skeleton of the perturbed triangulation
which is the complement of the handlebody neighborhood of the $n-k+1$ skeleton which avoids $K$.
(This proof was directly inspired by a conversation with Chris Bishop at Stony Brook University.)

\item[iii)] Now dual cells to a triangulation thicken to handles with two properties
\begin{itemize}
\item[a)] the handles have small diameters if the mesh of the triangulation is small.
\item[b)] the handles associated to the dual k-skeleton are convex and
only intersect k+1 at a time.  See Figure 1.
\end{itemize}

\item[iv)] Taking the dual handles that actually touch $K$ yields small neighborhoods of $K$  by iii a)
whose cohomology vanishes in dimensions above $k$ by iii b). These are
the $k$-good neighborhoods. 

\end{itemize}

\begin{figure}
\centerline{\includegraphics*[scale=0.5]{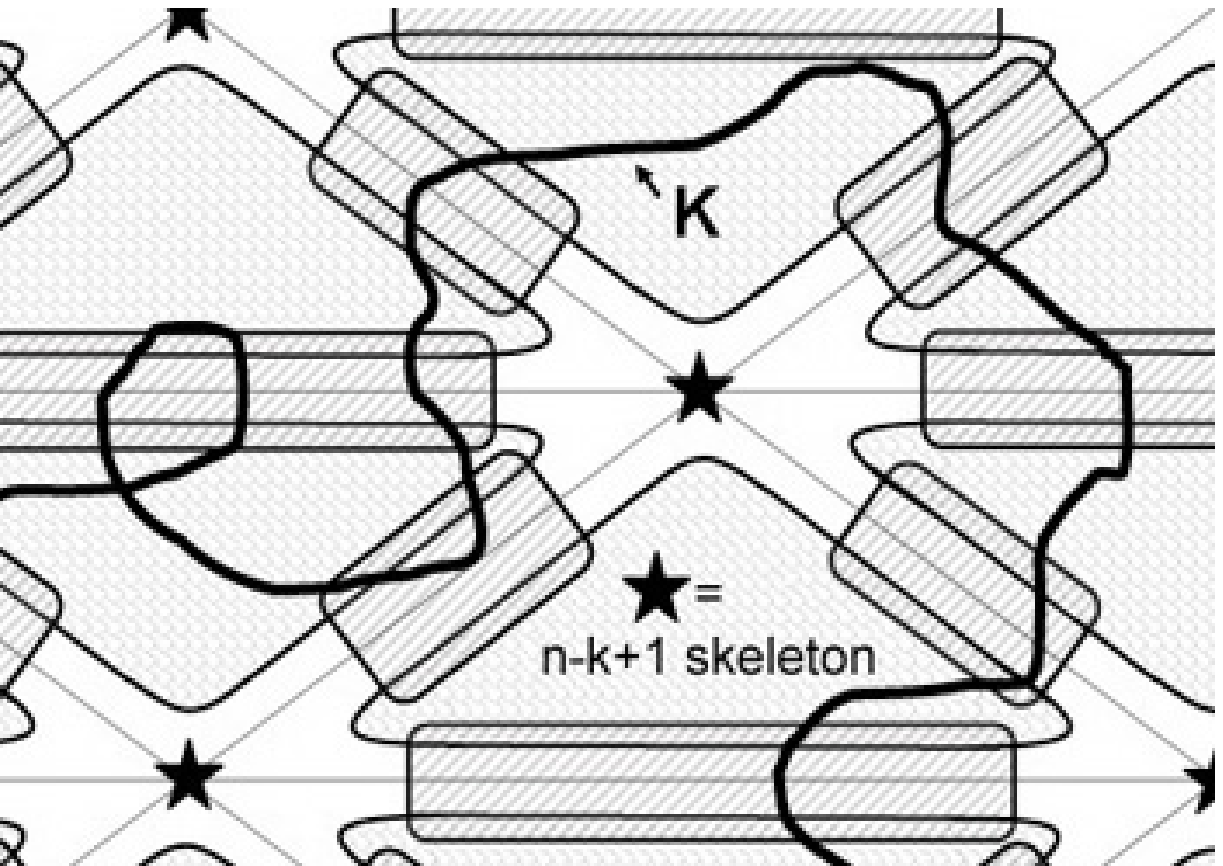}}
\caption{}
\end{figure}

\underline{\textbf{Fact 2.2}}:  If dim $M$ is at least $k$, a smooth
singular $(k-1)$ cycle is homologous in every neighborhood of its image to the
fundamental cycle of a piecewise smoothly embedded oriented $(k-1)$
pseudomanifold.  (See proof for a recall of the definition of pseudomanifold.)

\underline{\textbf{Remark}}:  Here we are considering (as did
Poincar\'{e} [1] in his original paper) smooth curvilinear nondegenerate
simplices in $M$ defined by smooth equations and smooth inequalities.
We consider $M$ as decomposed into such simplices, further arbitrarily
small diameter subdivisions, after Munkres [3], and perturbations of
these as in Figures 2 and 3.

\begin{figure}
\centerline{\includegraphics*[scale=0.5]{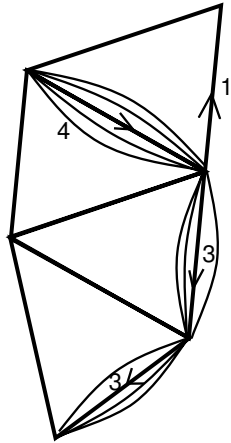}}
\caption{Splitting}
\end{figure}

\begin{figure}
\centerline{
\includegraphics*[scale=0.4]{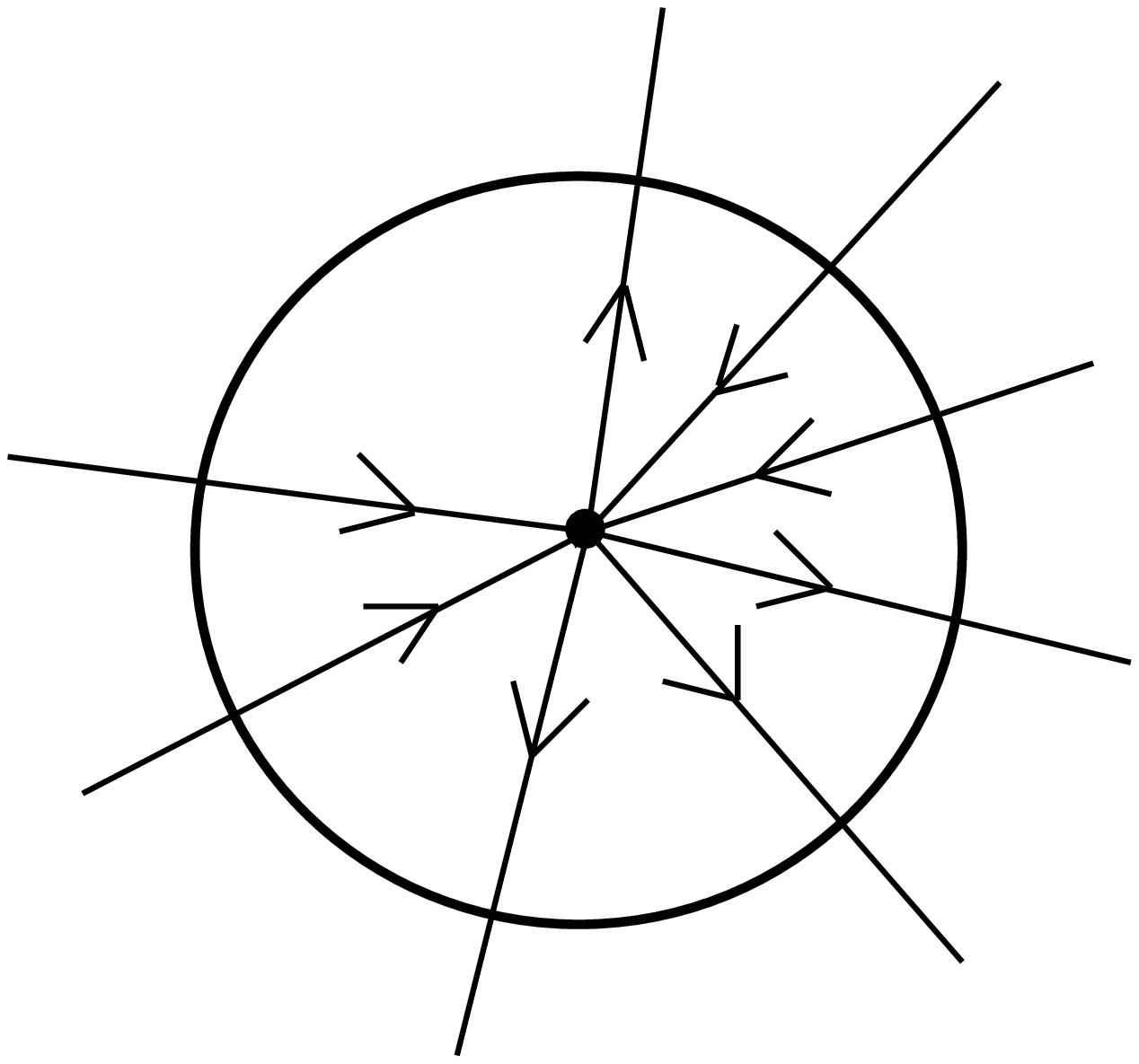}
\includegraphics*[scale=0.4]{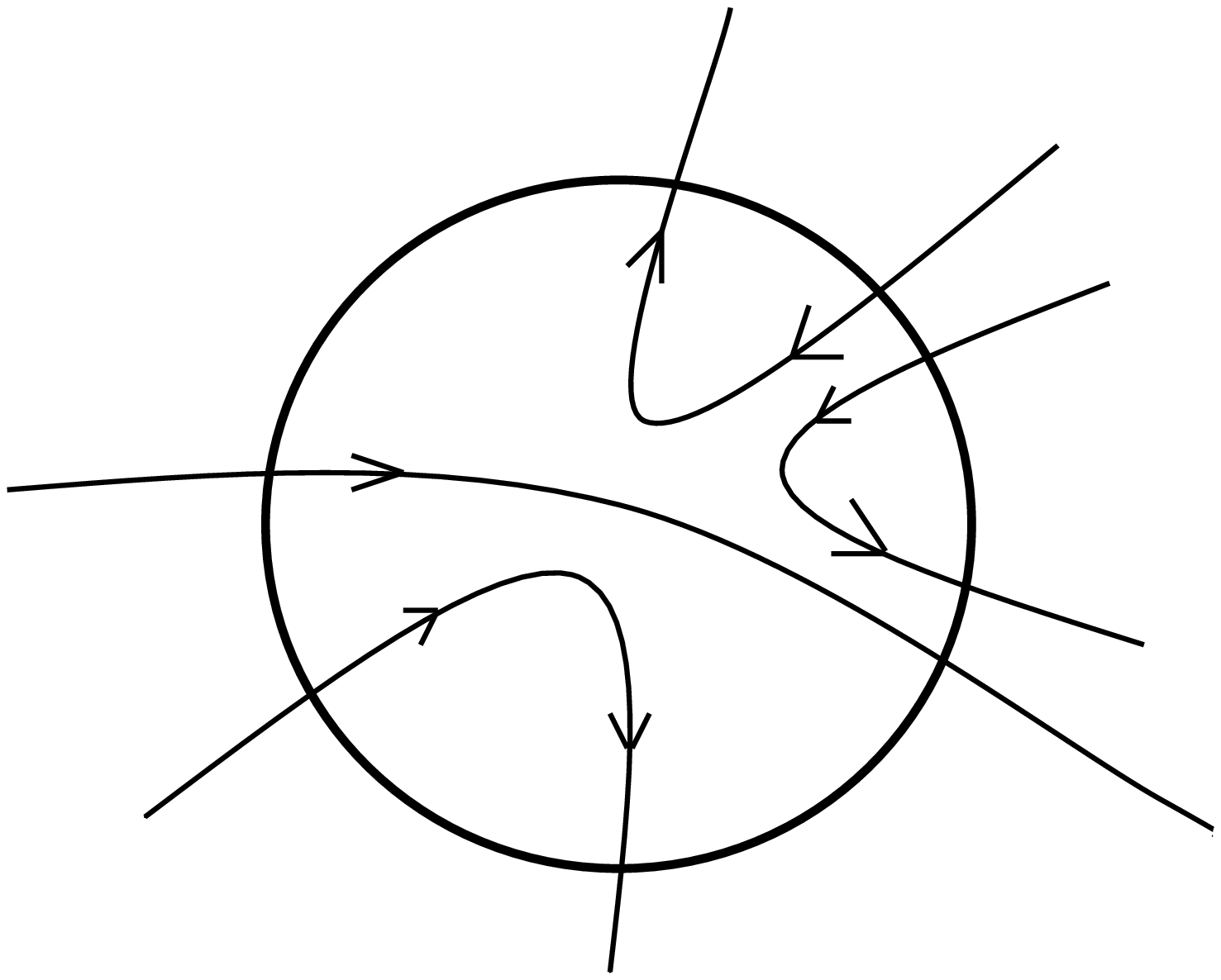}
}
\centerline{3a \hspace{5.5cm} 3b}
\caption{Resolving}
\end{figure}

If $k-1$ were the dimension of $M$, and $M$ were a closed oriented
manifold, the theorem would only be true for the generators of
homology.  Multiples of generators would not admit embedded representations.

\underline{\textbf{Proof}}:  
\begin{itemize}
\item[i)]  Let $W$ be a smooth manifold with boundary, chosen to be a
neighborhood of the image of our $(k-1)$ singular cycle inside the
given neighborhood.  We can choose a piecewise smooth triangulation of
$M$ in which $W$ is a subcomplex, cf. [6].

\item[ii)]  The smooth singular homology of $W$ is naturally
isomorphic (by the inverse of the inclusion) to the cellular homology
of $W$ for this triangulation.  Thus our cycle is homologous in $W$
(by a smooth singular homology) to a cellular $k-1$ cycle in $W$.

\item[iii)]  Now we orient the $(k-1)$ cells of the support of this
cycle so that the cycle is written with positive coefficients.

\item[iv)]  If a coefficient of a cell is an integer $n$ greater than
one we replace that cell by $n$ deformed copies, as in Figure 2.  This
uses the hypothesis $(k-1) < \textrm{dimension } M$.

\item[v)]  Now transverse to the $(k-2)$ cells of the deformed cycle
we have a picture like Figure 3a.  We now perform the deformation
indicated by Figure 3b.

\item[vi)]  The deformations of iv) and v) provide a further homology
of our cycle to one which is carried by the cell sum of a
$(k-1)$-dimension oriented curvilinear polyhedron so that each $(k-2)$
cell is the face of exactly two $(k-1)$ cells with cancelling
orientation.  Note:  This explains the phrase ``the fundamental cycle
of an oriented embedded pseudomanifold''.

\end{itemize}

\underline{\textbf{Fact 2.3}}:  Suppose dim $M \geq k$.  If the fundamental cycle $\xi$ of an
embedded oriented $(k-1)$ pseudomanifold is homologous to zero in $M$,
then $\xi$ is also bounds in some $(k-1)$-good neighborhood of its support.

\underline{\textbf{Proof}}:

\begin{itemize}
\item[i)]  Extend a cellular decomposition of the cycle to a piecewise
smooth cell decomposition of $M$.  Then $\xi = \partial w^{\prime
\prime}$ for some $k$-chain in this subdivision.

\item[ii)]  Orient the $k$ cells of $w^{\prime \prime}$ to get
positive integer coefficients.  If the dimension $M=k$ and $\xi$
bounds $w^{\prime\prime}$ in $M$, then $|\xi|$ separates $M$ and
$w^{\prime\prime}$ picks out a compact part of the
complement which is oriented and has oriented boundary equal to
$\xi$.  If the dimension $M \geq k+1$, apply the deformation iv) of the
previous proof, indicated in Figure 2, to get a new cellular homology
$w^{\prime}$ with all coefficients equal to 1.  We still have
$\partial w^{\prime} = \xi$.

\item[iii)]  Continuing in the case dimension $M \geq k+1$, apply the
deformation of v) of the previous proof to $(k-1)$ cells of
$|w^{\prime}|$ in the ``interior'' i.e. those not in $|\xi^{\prime}|$
on the boundary of $|w^{\prime}|$.  
At a boundary cell we perform a variant of this deformation at each
$(k-1)$ cell indicated by Figure 4, where the $(k-1)$ cell of $\xi$ is
depicted as a point.

\begin{figure}
\centerline{
\includegraphics*[scale=0.5]{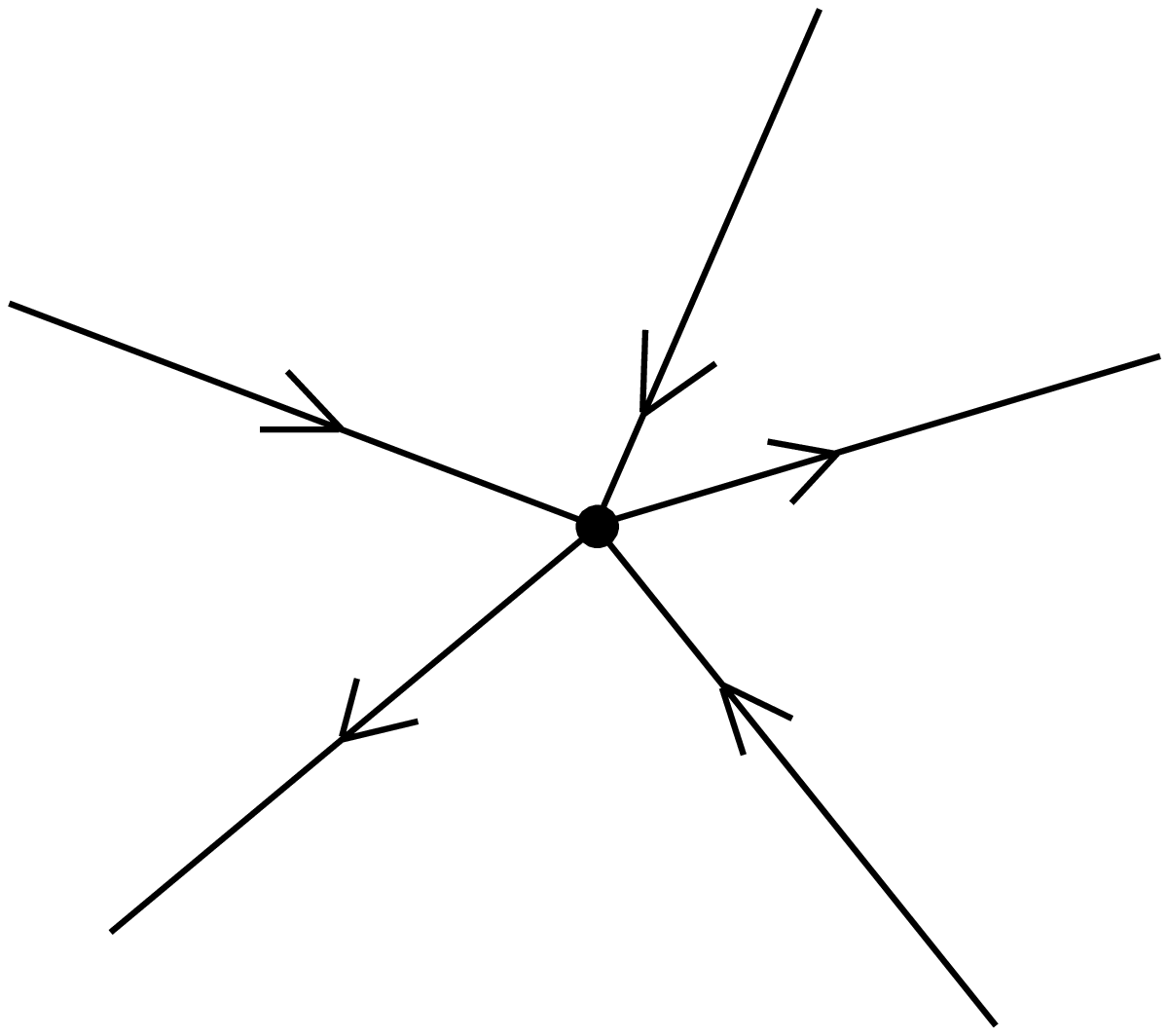}
\includegraphics*[scale=0.5]{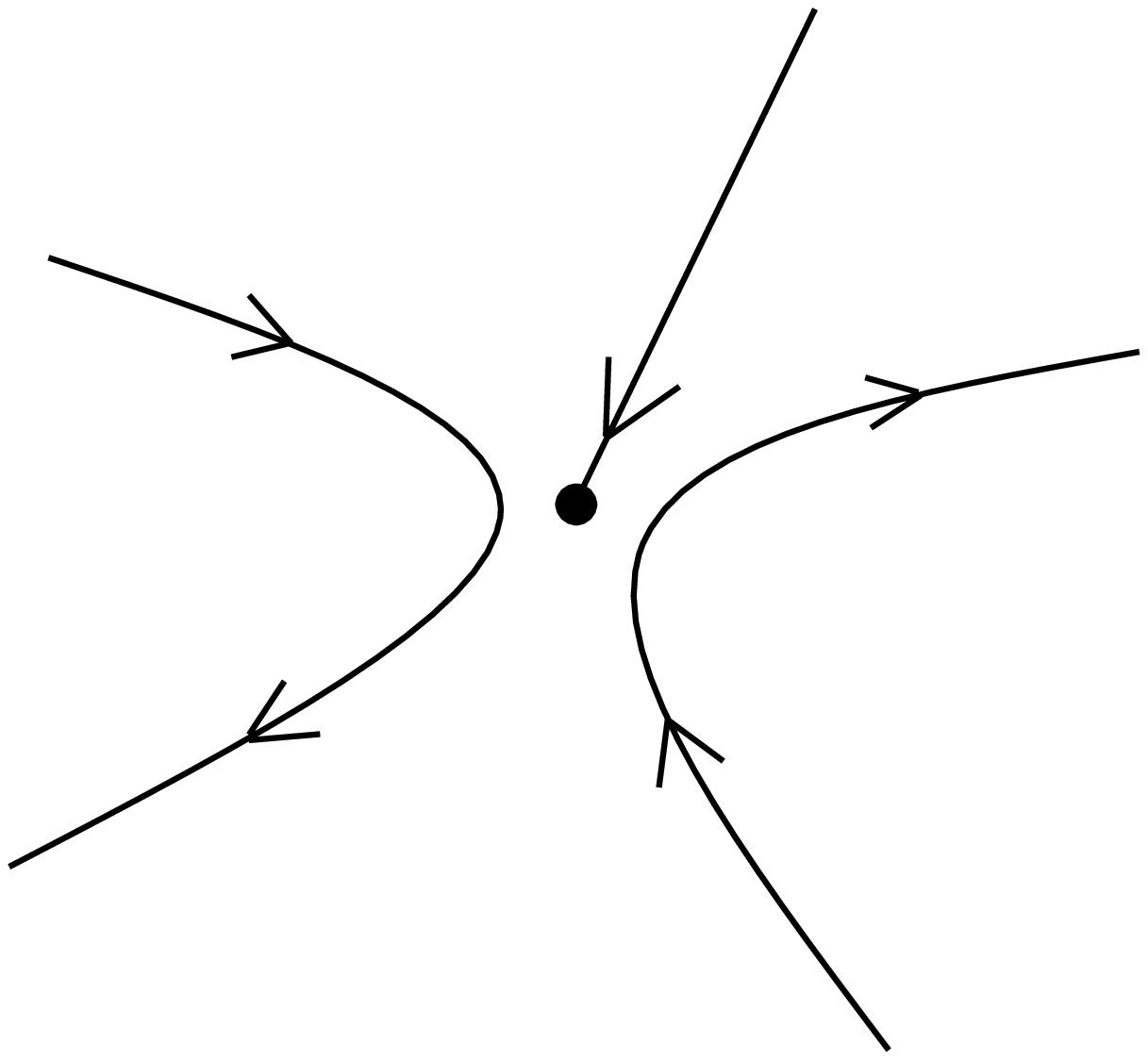}
}
\caption{}
\end{figure}

Conclusion:  In either case dimension $M = k$ or dimension $M \geq
k+1$ we obtain an embedded oriented pseudomanifold with boundary $\xi$.

\item[iv)]  Now the argument that a connected oriented triangulated
$k$-dimensional manifold with boundary deformation retracts to a
$(k-1)$ dimensional subcomplex goes by pushing in $k$-cells starting
from the boundary until none are left.  There are always exposed
$(k-1)$ faces to push in on until we get to a subcomplex of the $(k-1)$
skeleton because if not there would be a non-trivial $k$-cycle.  The
latter is impossible because we started with a connected oriented
manifold with boundary.  (The same argument works for unoriented
manifolds using mod 2 homology considerations.)

\item[v)]  The above arguments work word for word the same for
irreducible $k$-pseudomanifolds with boundary (i.e. the dual 1-skeleton
is connected) to show that such
an object deformation retracts to a $(k-1)$ dimensional subcomplex for
some curvilinear triangulation.

\item[vi)]  Now a regular neighborhood of the $k$-pseudomanifold with,
boundary has the same homotopy type and we have proved Fact 2.3.

\end{itemize}

\section*{References}

\begin{enumerate}
\item Poincar\'{e}, Henri.  ``Complement a l'Analysis Situs''.
Oeuvres Completes Tome VI.  Gauthier Villars Paris. 1954. p. 291.

\item Hurewicz, Witold and Wallman, Henry. ``Dimension
Theory''. Princeton University Press.  1941.

\item Munkres, J.R. ``Elementary Differential Topology''. Annals of
Mathematics Studies. Number 54. 1969. p. 79.

\item Sullivan, Dennis. ``Singularities in Spaces''.  Proceedings of the
Liverpool Conference on Singularities Symposium II, Liverpool, England.  Springer
Lecture Notes (209) pp. 196-206.

\item Chern, S.S. and Simons, James.  ``Characteristic Forms and Geometric Invariants''.
Annals of Mathematics. Vol. 99. No. 1. January 1974. pp. 48-69.

\item Simons, James.  ``Characteristic Forms and Transgression: Characters
Associated to a Connection''.  Stony Brook University Preprint.  1972.

\item Cheeger, Jeff. ``Multiplication of Differential Characters''.
Instituto Nazionale di Alta Mathematica, Symposia Mathematica
XI. 1973.  pp. 441-445.

\item Cheeger, Jeff and Simons, James.  ``Differential Characters and Geometric Invariants''.
Notes of Stanford Conference 1973, Lecture Notes in Math. No. 1167.
Springer-Verlag, New York.  1985.  pp. 50-90.

\item Gillet, H. and Soule, C..  ``Arithmetic Chow Groups and
Differential Characters, in Algebraic K-Theory: Connections with
Geometry and Topology''.  Jardine and Snaith, editors.  Kluwer Academic
Publishers.  1989.  pp. 30-68.

\item Harris, B.  ``Differential Characters and the Abel-Jacobi map, in
Algebraic K-Theory: Connections with Geometry and Topology''.  Jardine
and Snaith, editors.  Kluwer Academic Publishers.  1989. pp. 69-86.

\item Brylinski, Jean-Luc.  Loop Spaces, Characteristic Classes, and Geometric
Quantization.  Birkhauser, Boston. 1993.

\item Harvey, Reese; Lawson, Blaine; and Zweck, John.  ``The
de Rham-Federer theory of differential characters and character
duality''.  American Journal of Mathematics.  Vol 125.  2003.  pp. 791-847.

\item Harvey, Reese and Lawson, Blaine.  ``From Sparks to
Grundles-Differential Characters''. Communications in Analysis and
Geometry. Vol. 14, No. 1. 2006. pp. 25-58.

\item Hopkins, Michael and Singer, Isadore. ``Quadratic functions in
 geometry, topology,and M-theory''. http://arxiv.org/pdf/math/0211216.

\end{enumerate}

\end{document}